\documentstyle[11pt]{article}   

\newtheorem{th}{Theorem}   
\newtheorem{ax}{Axiom}  
\newtheorem{lm}{Lemma} 
\newtheorem{df}{Definition}    
\newtheorem{pr}{Proposition} 
\newtheorem{cl}{Corollary}  
\newtheorem{re}{Remark}    
\newtheorem{as}{Assumption}  
\newtheorem{ex}{Example}
\newcommand{\bth}{\begin{th}\hspace{-5pt}{\bf .} \ } 
\newcommand{\eth}{\end{th}}
\newcommand{\bax}{\begin{ax}\hspace{-5pt}{\bf .} \ } 
\newcommand{\eax}{\end{ax}}
\newcommand{\blm}{\begin{lm}\hspace{-5pt}{\bf .} \ }
\newcommand{\elm}{\end{lm}}
\newcommand{\bdf}{\begin{df}\hspace{-5pt}{\bf .} \ }   
\newcommand{\edf}{\end{df}} 
\newcommand{\bpr}{\begin{pr}\hspace{-5pt}{\bf .} \ } 
\newcommand{\epr}{\end{pr}}
\newcommand{\bcl}{\begin{cl}\hspace{-5pt}{\bf .} \ } 
\newcommand{\ecl}{\end{cl}}
\newcommand{\bre}{\begin{re}\hspace{-5pt}{\bf .} \ }
\newcommand{\ere}{\end{re}}
\newcommand{\bas}{\begin{as}\hspace{-5pt}{\bf .} \ }
\newcommand{\eas}{\end{as}}        
\newcommand{\bex}{\begin{ex}\hspace{-5pt}{\bf .} \ }
\newcommand{\eex}{\end{ex}}
\newcommand{\bpf}{\noindent {\it Proof}\hspace{0.1truecm}: }
\newcommand{\epf}{\hfill${\scriptstyle\diamondsuit\diamondsuit
\diamondsuit}$\par\vspace{1.8mm}\noindent}
\newcommand{\bit}{\begin{itemize}}
\newcommand{\eit}{\end{itemize}\par\noindent}
\newcommand{\beq}{\begin{equation}} 
\newcommand{\eeq}{\end{equation}\par\noindent}
\newcommand{\beqa}{\begin{eqnarray*}}
\newcommand{\eeqa}{\end{eqnarray*}\par\noindent}
\newcommand{\beqn}{\begin{eqnarray}}
\newcommand{\eeqn}{\end{eqnarray}\par\noindent}
\renewcommand{\a}{{\cal A}}
\renewcommand{\b}{{\cal B}}
\renewcommand{\c}{{\cal C}}

\newcommand{\e}{{\cal E}}

\renewcommand{\j}{{\cal J}}

\renewcommand{\l}{{\cal L}}

\newcommand{\p}{{\cal P}}
\newcommand{\q}{{\cal Q}}

\newcommand{\s}{{\cal S}}

\newcommand{\1}{{\bf 1}}

\title{{\bf State Transitions as Morphisms \\ for Complete Lattices}} 
\date{}
\author{}

\begin{document}  
\pagenumbering{arabic}      
\pagestyle{plain}
\maketitle 

\par\vspace{-2cm}\par
\begin{center}
{\it Bob Coecke}\\
{\scriptsize FUND-DWIS, Free University of Brussels,\\
Pleinlaan 2, B-1050 Brussels, Belgium.\\
bocoecke@vub.ac.be}
\par\medskip\par
{\it Isar Stubbe}\\
{\scriptsize AGEL-MAPA, Universit\'e Catholique de Louvain,\\
Ch$.$ du Cyclotron 2, B-1348 Louvain-La-Neuve, Belgium.\\
i.stubbe@agel.ucl.ac.be}   
\par\bigskip\par
Appeared in {\it International Journal of Theoretical physics}\ \ {\bf
39(3)}, 605--614 (2000).  
\end{center}    

\begin{abstract}  We enlarge the hom-sets of categories of
complete lattices by introducing `state transitions' as
generalized morphisms.  The obtained category will then be
compared with a functorial quantaloidal enrichment and a
contextual quantaloidal enrichment that uses a specific concretization in
the category of sets and partially defined maps ($\p arset$).
\end{abstract}

\section{Introduction}

In this paper we present a construction that abstracts the concept of
`state transition' as introduced in (Amira {\it et al.}, 1998) and
(Coecke and Stubbe, 1999a; 1999b), making it applicable to any
subcategory $\a$ of ${\cal JCL}at$. We compare this construction, the
result of which is a quantaloid that we call $\q^{st}\a$, with two other
quantaloidal extensions that arise naturally when considering the action
of the power-functor on
$\a$. In fact, one of these natural extentions is functorial, we denote
it by
$\q^{-}\a$, and the other, called $\q^{+}\a$, is contextual in the sense
that its construction relies entirely on the fact that $\a$ is a $\p
arset$-concrete category. The main result of this paper is then that in
all non-trivial cases $\q^{st}\a$ lies strictly between $\q^{-}\a$ and
$\q^{+}\a$.

Applying the construction $\q^{st}$ to the category $\p rop$,
which was introduced in (Moore, 1995; 1999), reveals that the latter has
to be `enriched' in order to constitute an appropriate mathematical
object for defining `state transitions'. However, it must be noted that
the physical inspiration for `our' categories is essentially different
from Moore's: $\p rop$ is constructed to express the equivalence of a
property lattice and a state space description for a
physical system as a categorical equivalence of $\p rop$ and $\s tate$,
inspired by a similar situation for categories for projective geometries
(Faure and Fr\"olicher, 1993; 1994), and delivers the
mathematical context that embeds the proof of Piron's representation
theorem (Piron, 1964; 1976). On the contrary, we propose a
dualization/generalization for the notion of property transition,
previously motivated to be a join-preserving map (Pool, 1968; Faure {\it
et al.}, 1995; Amira {\it et al.}, 1998).
  
\section{`Enriching' $\p arset$-concrete categories}  

For a general overview of the theory of categories we refer to (Adamek
{\it et al.}, 1990; Borceux, 1994).  For quantaloids we refer to (Pitts,
1988; Rosenthal, 1991).
\bdf A `quantaloid' is a category such that: 
\par\noindent (i) every hom-set is a complete lattice, its join is
usually denoted by
$\bigvee$;
\par\noindent (ii) composition of morphisms distributes on both sides
over arbitrary joins.
\par\noindent  Let $\q $ and $\q '$ be quantaloids. A `quantaloid
morphism' from $\q $ to
$\q '$ is a functor
$F:\q \rightarrow\q '$ such that on hom-sets it induces join-preserving
maps
$\q (A,B)\rightarrow\q '(FA,FB)$.
\edf   For example, $\j\c\l at$ is the quantaloid of complete lattices
and join-preserving maps in which the join of maps is computed
pointwise. A quantaloid with one object is commonly known as a `unital
quantale', a quantaloid morphism between two one-object quantaloids is
known as a `unital quantale morphism' (Rosenthal, 1990). Any subcategory
of a quantaloid that is closed under the inherited join of morphisms, is a
subquantaloid. Thus any full subcategory of a quantaloid is a 
subquantaloid, and selecting from a given quantaloid certain morphisms
but keeping all the objects, gives rise to a subquantaloid if and only
if the inherited join of morphisms is  internal.
 
Let $\p arset$ denote the category of sets $A,B,...$ and partially
defined functions
$f:A\setminus K\to B$ --- where $K\subseteq A$ is  called the `kernel'
of the partially defined function $f$, also written as ${\sf ker}f$. Then
the power-functor is defined as
\[  
\begin{array}{ccl}
\p & : & \p arset\to\j\c\l at\\
 & : &
\left\{\begin{array}{rcl}   A & \mapsto &\ 2^A \\ f:A\setminus K\to B &
\mapsto & \p f:2^A\to 2^B:T\mapsto\{f(t)\mid t\in T\setminus K\}
\end{array}
\right.
\end{array}
\]  that is, sets are mapped onto their powersets, partially defined
functions are mapped onto the `direct image mapping' which is indeed a
union-preserving map. This functor is faithful and injective, but neither
full nor surjective.

For any $\p arset$-concrete category $\a$, that is a category that comes
equipped with a faithful functor $U:\a\to\p arset$, we can compose
functors
$\a\stackrel{U}{\longrightarrow}\p
arset\stackrel{\p}{\longrightarrow}\j\c\l at$ and use this to define a
category that we shall denote by $\p\a$: it has the same objects as $\a$,
and for the hom-sets we define:
$$\p\a(A,B)=\{(\p\circ U)(f)\mid f\in \a(A,B)\}.$$ So in particular the
hom-sets are posets for the pointwise order. Any functor $F:\a\to \b$
between two such categories
$\a\stackrel{U}{\to}\p arset$ and $\b\stackrel{V}{\to}\p arset$ defines a
functor
\[  
\begin{array}{ccl}
\p F & : & \p\a\to\p\b\\
 & : &
\left\{\begin{array}{rcl}   A & \mapsto & FA \\
\p Uf:2^{UA}\to 2^{UB} & \mapsto & \p VFf:2^{VFA}\to 2^{VFB}
\end{array}
\right.
\end{array}
\] which goes to say that $\p$ is functorial on the quasi-category of $\p
arset$-concrete categories\footnote{Remark that $\s et$ is the
subcategory of $\p arset$  with the same objects but of which all
morphisms have an empty kernel, and that the domain restriction of $\p$
to $\s et$ yields that, for any morphism
$f:A\to B$ in $\s et$, $\p f:2^{UA}\to 2^{UB}$ is a union-preserving map
that maps only the emptyset on the emptyset: $(\p
f)(T)=\emptyset\Leftrightarrow T=\emptyset$.}.

Since $\p\a(A,B)\subseteq\j\c\l at(2^{UA},2^{UB})$ for any two
$\p\a$-objects
$A$ and $B$, we can define a category $\q^+\a$ with the same objects as
$\p\a$ (thus the same objects as $\a$) but of which the hom-sets are
precisely $$\q^+\a(A,B)=\j\c\l at(2^{UA},2^{UB}).$$ Explicitly this means
that a morphism
$f:A\to B$ in $\q^+\a$ is determined by an underlying union-preserving
map
$f:2^{UA}\to 2^{UB}$. So $\q^+\a$ is a quantaloid with respect to the
pointwise union of maps. Secondly, by $\q^-\a$ we shall denote the
category with the same objects as $\p\a$ (thus the same objects as $\a$
and as $\q^+\a$), but of which the hom-sets are the complete lattices
that one obtains if one closes the
$\p\a$-hom-sets for (arbitrary) pointwise unions of maps: a morphism 
$f:A\to B$ in $\q^-\a$ is thus determined by an underlying map\,
$\bigcup_i\p Uf_i$,  the join of maps being their pointwise union,
for a set of 
$\a$-morphisms $\{f_i:A\to B\}_i$. If $\a$ is a category in which every
hom-set contains a non-zero element --- that is, for any two objects
$A$, $B$ of $\a$ there is an $f\in\a(A,B)$ such that
${\sf ker}Uf\neq {\sf dom}Uf$ --- then any hom-set of $\p\a$ contains a
non-zero map, thus any hom-set of
$\q^-\a$ is a complete lattice (the condition on the hom-sets of $\a$
makes sure that any hom-set of $\q^-\a$ contains at least distinct bottom
and top which we require for any complete lattice), so
$\q^-\a$ is a quantaloid with respect to pointwise union of morphisms.
For a functor
$F:\a\to\b$ between two such categories
$\a\stackrel{U}{\to}\p arset$ and $\b\stackrel{V}{\to}\p arset$, we can
define a corresponding quantaloid morphism
\[
\begin{array}{ccl}
\q^-F & : & \q^-\a\to\q^-\b\\
 & : &
\left\{\begin{array}{rcl}   A & \mapsto & FA \\
\bigcup_i\p Uf_i:2^{UA}\to 2^{UB} & \mapsto & \bigcup_i\p
VFf_i:2^{VFA}\to 2^{VFB}
\end{array}
\right.
\end{array}
\] which translates the idea that the construction of $\q^-\a$ is
functorial. Note that no such functor can exist for the $\q^+$-case
because in $\q^+\a$ there may be morphisms for which there exists no
connection at all with $\a$-morphisms! So we could say that $\q^-\a$ is
a `functorial' enrichment of $\a$, and $\q^+\a$ merely a `contextual'
enrichment.

By construction, $\q^-\a$ is a subquantaloid of
$\q^+\a$ since the hom-sets of the former are join-sublattices of the
hom-sets of the latter (they have the same objects). For an important
class of categories this inclusion is strict: we will consider a
category $\a$ of which the objects are bounded posets (partially ordered
sets with a greatest element
$1$ and a least element $0$, these elements being different) and the
morphisms are isotone mappings that map least elements onto least
elements (that is, $0\mapsto 0$), together with a forgetful functor
\[
\begin{array}{ccl}  
U & : & \a\to\p arset\\
 & : & \left\{\begin{array}{rcl}  
(A,\leq) & \mapsto & A_0 \\   f:(A,\leq)\to(B,\leq) & \mapsto &
Uf:A_0\setminus K\to B_0:t\mapsto f(t)
\end{array}
\right.
\end{array}
\]    where $K=\{t\in A_0\mid f(t)=0\}={\sf ker}Uf$ and
$A_0=A\setminus\{0\}$. From now on, we will always assume that every
hom-set has at least one non-zero element, this to assure that the
construction
$\q^-\a$ yields a quantaloid.
\bpr\label{Q^-versusQ^+} For such a category $\a$,
$\q^-\a=\q^+\a$ if and only if all objects of $\a$ are 2-chains.
\epr
\bpf
If $\a$ contains an object $(A,\leq)$ in which there exists
$0<a<1$, then
\[
\begin{array}{cc} f:2^{A_0}\to 2^{A_0}: &
\left\{\begin{array}{lcr}   T\mapsto \{a\} & \Leftrightarrow & a\in T \\
T\mapsto \emptyset & \Leftrightarrow & a\not\in T \\
\end{array}
\right.
\end{array}  
\] is trivially a union-preserving map, thus $f\in \q^+\a(A,A)$. For
this $f$ to be a morphism of $\q^-\a$ as well, there must be a set of
$\a$-morphisms $h_i:A\to A$ such that precisely $f=\bigcup_i\p U h_i$:
$$\forall T\in 2^{A_0}:f(T)=\cup_i\{Uh_i(t)\mid t\in T\setminus {\sf ker}
Uh_i\}.$$ Since $f(\{1\})=\emptyset$, we must have that $h_i(1)=0$ for all
those $h_i$, but on the other hand $f(\{a\})=\{a\}$, which
means that at least one
$h_k$ is such that $h_k(a)=a$. But this contradicts with the
isotonicity of the $\a$-morphism $h_k$.
If on the contrary the category $\a$ contains only posets which are
(isomorphic to) 2-chains then it is merely an observation that the
functorial and the contextual enrichment are the same.
\epf 
As an example, consider:
\[  
\begin{array}{ccl} 
U & : & \j\c\l at\to\p arset\\
 & : & \left\{\begin{array}{rcl}  
(L,\vee) & \mapsto & L_0 \\   f:(L,\vee)\to(M,\vee) & \mapsto &
Uf:L_0\setminus {\sf ker}Uf\to M_0:t\mapsto f(t)
\end{array}
\right.
\end{array}
\]
With respect to this faithful functor we have the following categories, all
three with as objects complete lattices:
\bit
\item
$\p\j\c\l at$: a morphism $f:(L,\vee)\to (M,\vee)$ corresponds to the image
$\p Uf:2^{L_0}\to 2^{M_0}$ of a join-preserving map $f:L\to M$ (thus the
underlying map preserves unions).
\item
$\q^-\j\c\l at$: a morphism $f:(L,\vee)\to (M,\vee)$ corresponds to a map
$f:2^{L_0}\to 2^{M_0}$ that can be written as the pointwise union of
$\p\j\c\l at$-maps (then it automatically preserves unions).
\item
$\q^+\j\c\l at$: a morphism $f:(L,\vee)\to (M,\vee)$ corresponds to a map
$f:2^{L_0}\to 2^{M_0}$ that preserves unions.
\eit As a corollary of proposition \ref{Q^-versusQ^+} the
quantaloid inclusion $$\q^-\j\c\l at\hookrightarrow \q^+\j\c\l at$$ is
strict. In the next section we introduce the notion of `state
transition' as morphism between complete lattices, and show that these
can organize themselves in a quantaloid that lies between $\q^-\j\c\l
at$ and $\q^+\j\c\l at$.

\section{`State transitions' as morphisms between complete lattices}
 
Given two complete lattices $(L,\vee)$, $(M,\vee)$ and a
map $f:2^{L_0}\to 2^{M_0}$ that preserves unions, it is easy to verify that
if there exists a map $g:L\to M$ that makes 
\[
\begin{array}{ccc} L & \stackrel{g}{\longrightarrow} & M \\
\uparrow\scriptstyle{\vee} & & \uparrow\scriptstyle{\vee} \\ 2^{L_0} &
\stackrel{f}{\longrightarrow} & 2^{M_0} 
\end{array}
\] commute --- the vertical uparrows denoting, for the respective
lattices, the maps
$T\mapsto\vee T$ --- then this map $g$ is a unique $\j\c\l at$-morphism.
Therefore the following definition makes sense:
\bdf\label{statetransition} Given a subcategory $\a$ of $\j\c\l at$ and
two of its objects $(L,\vee)$,
$(M,\vee)$, a union-preserving map
$f:2^{L_0}\to 2^{M_0}$ is called `state transition with respect to $\a$'
if there exists an
$\a$-morphism $f_{pr}:L\to M$ that makes the diagram 
\[
\begin{array}{ccc} L & \stackrel{f_{pr}}{\longrightarrow} & M \\
\uparrow\scriptstyle{\vee} & & \uparrow\scriptstyle{\vee} \\ 2^{L_0} &
\stackrel{f}{\longrightarrow} & 2^{M_0} 
\end{array}
\] commute. This (unique) morphism
$f_{pr}$ is then called the `property transition corresponding to $f$'.
\edf  For a given subcategory $\a$ of $\j\c\l at$ all morphisms of
$\p\a$ are state transitions with respect to $\a$:
\bpr\label{propstatetransition} For every $\a$-morphism $f$ we have that
$(\p Uf)_{pr}=f$. For every set $\{f_i\}_{i\in I}$ of parallel
$\a$-morphisms we have that $(\bigcup_i\p Uf_i)_{pr}=\bigvee_if_i$, these
joins being computed pointwise.
\epr
\bpf For $f\!:\!(L,\vee)\!\to\!(M,\vee)$ in $\a$ the direct image map $\p
Uf:2^{L_0}\to2^{M_0}$ is union-preserving, and for any $T\in 2^{L_0}$:
$\vee\p Uf(T)  =  \vee\{Uf(t)\mid t\in T\setminus {\sf ker}Uf\}
           =  \vee\{f(t)\mid t\in T\}
           =  f(\vee T)$.
The second assumption can be proven likewise.
\epf   This proposition does not say that any $\q^-\a$-morphism ---
which is by definition an arbitrary pointwise union of $\p\a$-morphisms
--- is a state transition with respect to $\a$, for it may very well be
that, even with all $f_i$ in $\a$, the corresponding property transition
$(\bigcup_i\p Uf_i)_{pr}=\bigvee_if_i$ is not an
$\a$-morphism. So in general any $\q^-\a$-morphism is a state transition
only with respect to $\j\c\l at$. But it follows immediately that:
\bcl\label{quantforMINCROSS} Every $\q^-\a$-morphism is a state
transition with respect to $\a$ if and only if $\a$ is a subquantaloid
of ${\cal JCL}at$.
\ecl For any state transition $f$ (the subcategory $\a$ is of no
importance here) and any $T,T'\in 2^{L_0}$ we have
\beq\label{respectsresolution}
\vee T=\vee T'\Rightarrow \vee f(T)=\vee f(T'), 
\eeq because $\vee f(T)=g(\vee T)=g(\vee T')=\vee f(T')$. But also a
converse is true: any   union-preserving map $f:2^{L_0}\to 2^{M_0}$ that
meets eq(\ref{respectsresolution}) is a state transition with respect to
$\j\c\l at$, since in this case the map
$$f_{pr}:L\to M:t=\vee T\mapsto \vee f(T)$$ is a well-defined $\j\c\l
at$-morphism that makes the square commute. Thus:
\bpr\label{previousdefinition} A map $f:2^{L_0}\to 2^{M_0}$ is a state
transition with respect to $\j\c\l at$ if and only if it preserves
unions and meets the condition of eq(\ref{respectsresolution}).
\epr In fact, this proposition was taken as definition for `state
transition (with respect to $\j\c\l at$)' in (Coecke and Stubbe, 1999a;
1999b) and, under a slightly different form, in (Amira {\it et al.},
1998). Remark that the conditions in this proposition are also necessary
for $f$ to be a state transition with respect to
$\a$, a subcategory of $\j\c\l at$, but in general not sufficient. This
will prove its use in proposition \ref{kuuuplus}.

By pasting together commuting diagrams it is easily seen that the
composition of any two state transitions (with respect to a certain
$\a$) is again a state transitions (with respect to that same
$\a$) with corresponding property transition $(g\circ
f)_{pr}=g_{pr}\circ f_{pr}$. Trivially also the diagram with identities
commutes, thus identities on powersets are state transitions (with
respect to whatever $\a$) with corresponding property transition
$(id_{2^{L_0}})_{pr}=id_L$. 
\begin{center}
\[
\begin{array}{ccccc}  L & \stackrel{f_{pr}}{\longrightarrow} & M &
\stackrel{g_{pr}}{\longrightarrow} & N \\
\uparrow\scriptstyle{\vee} & & \uparrow\scriptstyle{\vee} & &
\uparrow\scriptstyle{\vee} \\  2^{L_0} & \stackrel{f}{\longrightarrow} &
2^{M_0} & \stackrel{g}{\longrightarrow} & 2^{N_0} 
\end{array}
\quad\quad
\quad\quad
\quad\quad
\begin{array}{ccc} L & \stackrel{id_L}{\longrightarrow} & L \\
\uparrow\scriptstyle{\vee} & & \uparrow\scriptstyle{\vee} \\ 2^{L_0} &
\stackrel{id_{2^{L_0}}}{\longrightarrow} & 2^{L_0} 
\end{array}
\] 
\end{center}  This means that for a given subcategory $\a$ of $\j\c\l
at$ we can define a category $\q^{st}\a$ with the same objects as $\a$,
in which a morphism $f:(L,\vee)\to(M,\vee)$ is determined by an
underlying state transition with respect to $\a$. This category contains
$\p\a$ (by proposition \ref{propstatetransition}) and is in turn a
subcategory of $\q^+\a$:
$$\p\a\hookrightarrow\q^{st}\a\hookrightarrow \q^+\a.$$ 
The foregoing is best summarized by the functor
\[  
\begin{array}{cc} F_{pr}:\q^{st}\a\to\a: & \left\{\begin{array}{rcl}  
(L,\vee) & \mapsto & (L,\vee) \\   f & \mapsto & f_{pr}
\end{array}
\right.
\end{array}
\] which is full because
$F_{pr}\circ\p:\a\to\a$ is the identity (see proposition
\ref{propstatetransition}); this fact expresses explicitly the
duality of state transitions and property transitions (with respect to
$\a$).

\section{State transitions and quantaloids}

If moreover $\a$ is a subquantaloid of $\j\c\l at$ then  commutation of
the left diagram for each $i\in I$ implies commutation of the right
diagram:
\begin{center}
\[
\begin{array}{ccc} L & \stackrel{f_{i,pr}}{\longrightarrow} & M \\
\uparrow\scriptstyle{\vee} & & \uparrow\scriptstyle{\vee} \\ 2^{L_0} &
\stackrel{f_i}{\longrightarrow} & 2^{M_0} 
\end{array}
\quad\quad
\quad\quad
\quad\quad
\begin{array}{ccc} L & \stackrel{\bigvee_if_{i,pr}}{\longrightarrow} & M
\\
\uparrow\scriptstyle{\vee} & & \uparrow\scriptstyle{\vee} \\ 2^{L_0} &
\stackrel{\bigcup_if_i}{\longrightarrow} & 2^{M_0} 
\end{array}
\]  
\end{center} since these joins of maps are computed pointwise, which
goes to say that, if $\a$ is a subquantaloid of $\j\c\l at$ then
$\q^{st}\a$ is a quantaloid and  the functor
$F_{pr}:\q^{st}\a\to\a$ is a full quantaloid morphism.
$\q^{st}\a$ is then a subquantaloid of $\q^+\a$, and it contains
$\q^-\a$ as subquantaloid (by virtue of corollary \ref{quantforMINCROSS}
and the fact that in all three quantaloids joins of morphisms are
computed as pointwise unions):
$$\q^-\a\hookrightarrow\q^{st}\a\hookrightarrow\q^+\a.$$ Next we point
out in detail the relations between these three different constructions.
\bpr\label{kuuuplus} For a subcategory $\a$ of $\j\c\l at$,
$\q^{st}\a=\q^+\a$ if and only if all objects of ${\cal A}$ are
2-chains.
\epr
\bpf If $\a$ contains an object in which there
is a chain $0<a<1$.
With the ``same'' counterexample as in the proof of proposition
\ref{Q^-versusQ^+} it can be seen that
$\vee f(\{1\})=\vee\emptyset=0\neq a=\vee\{a\}=\vee f(\{a,1\})$ although
$\vee\{1\}=\vee\{a,1\}$, which contradicts proposition
\ref{previousdefinition}, such that $f$ cannot be a state transition with
respect to $\j\c\l at$ and {\it a fortiori} $f$ cannot be a state
transition with respect to $\a$, a subcategory of $\j\c\l at$. Conversely,
if ${\cal A}$ contains only two-chains, then it is trivial that
$\q^{st}\a=\q^+\a$.
\epf
\bpr 
For a subcategory $\a$ of $\j\c\l at$,
$\q^-\a\not\supseteq\q^{st}\a$ if there is an $\a$-object
$(L,\vee)$ in which there exist elements $a,b,c$ such that $a<b\vee c$,
$a\not\leq b$, $a\not\leq c$.
\epr
\bpf Consider the map
\[
\begin{array}{cc} f:2^{L_0}\to 2^{L_0}: &
\left\{\begin{array}{l}  
\{t\}\mapsto\{t\} \ \ {\rm for}\ \ t\neq b\vee c \\
\{b\vee c\}\mapsto \{a,b,c\} \\ T\mapsto \cup_{t\in T}f(\{t\}) 
\end{array}
\right.
\end{array}  
\] for which clearly $f_{pr}=id_L$, thus it is a state transition (with
respect to whatever
$\a$). Let $g:L\to L$ be a join-preserving map that maps $a\mapsto a$,
then
$\p g\leq f\Rightarrow \p g(\{b\vee c\})\subseteq f(\{g\vee c\}) 
           \Rightarrow g(b)\vee g(c)\in\{a,b,c\} 
           \Rightarrow g(a)\leq g(b\vee c)=g(b)\vee g(c)\in\{a,b,c\}$
and   
$g(a)=a \Rightarrow a=g(b)\vee g(c)
        \Rightarrow g(b)\leq a \ \ {\rm and}\ \ g(c)\leq a$.  Should
$g(b)=a$ then $\p g(\{b\})=\{a\}$ and thus $\p g\not\leq f$. If on the
contrary $g(b)<a$ then $g(c)=a$, because
$g(c)<a$ would imply that $g(b)\vee g(c)<a$, but then $\p
g(\{c\})=\{a\}$ and thus $\p g\not\leq f$. Therefore, for any
$g\in\j\c\l at(L,L)$ such that $\p g\leq f$ we necessarily have that
$g(a)=0$. But then it is impossible to ever write a pointwise union
$f=\bigvee_i\p g_i$ with
$\{g_i\}_i\subseteq\j\c\l at(L,L)$, because the join on the right hand
side always fails to map
$\{a\}\mapsto\{a\}$.
\epf 
\bcl If $\a$ is a subquantaloid of $\j\c\l at$ that contains an object
in which there exist elements
$a,b,c$ such that $a<b\vee c$, $a\not\leq b$, $a\not\leq c$, then both
the inclusions of quantaloids
$\q^-\a\hookrightarrow\q^{st}\a\hookrightarrow\q^+\a$ are strict.
\ecl As $\j\c\l at$ is a trivial subquantaloid of $\j\c\l at$, we have
the full quantaloid morphism
\[  
\begin{array}{cc} F_{pr}:\q^{st}\j\c\l at\to\j\c\l at: &
\left\{\begin{array}{rcl}   (L,\vee) & \mapsto & (L,\vee) \\   f &
\mapsto & f_{pr}
\end{array}
\right.
\end{array}    
\]   expressing explicitly the duality between state transitions and
property transitions. Since
$\q^{st}\j\c\l at$  contains $\q^-\j\c\l at$, which in turn contains
$\q^-\a$ for any subcategory $\a$ of $\j\c\l at$, we can consider
$F_{pr}\q^-\a$, the image of $\q^-\a$ through this functor: this is the
smallest subquantaloid of $\j\c\l at$ that contais $\a$, it emerges by
closing all hom-sets
$\a(A,B)\subseteq\j\c\l at(A,B)$ for arbitrary (pointwise) joins.
Evidently, if $\a$ is a subquantaloid of $\j\c\l at$, then and only then
$\a=F_{pr}\q^-\a$. It can be verified straightforwardly that the
assignment $\a\mapsto F_{pr}\q^-\a$ is functorial; we will denote the
corresponding functor (that thus acts on subcategories of $\j\c\l at$
and functors between these) as
$\e$ and we will refer to
$\e\a$ as the {\it `pre-enrichment'} of
$\a$. Obviously the quantaloid $\q^-\a$ is included in $\q^-\e\a$, so we
can write the following inclusion of quantaloids as generalization of
the previous material:
$$\q^-\a\hookrightarrow\q^-\e\a\hookrightarrow\q^{st}\e\a\hookrightarrow\q^+\e\a.$$
Using the various previous propositions, conditions may be given for
these inclusions to be strict.

\section{Conclusion and examples}  

\bth      For any subcategory $\a$ of $\j\c\l at$ that contains an
object in which there exist elements $a,b,c$ such that $a<b\vee c$,
$a\not\leq b$, $a\not\leq c$, the quantaloid inclusions
$\q^-\e\a\hookrightarrow\q^{st}\e\a\hookrightarrow\q^+\e\a$ are strict.
Here $\e$ stands for the minimal extension of $\a$ to a quantaloid,
$\a\mapsto F_{pr}\q^-\a$. If moreover $\a$ is a subquantaloid of $\j\c\l
at$, then $\e\a=\a$.
\eth The essence of having an inclusion $\q^-\a\hookrightarrow\q^{st}\a$
--- or forcing it as
$\q^-\e\a\hookrightarrow\q^{st}\e\a$ --- should be understood in the
following way: The join
of maps in $\q^{st}\a$ physically stands for a lack of
knowledge on possible state transitions (Amira {\it et al.}, 1998; Coecke
and Stubbe, 1999a). Therefore, any general collection of state
transitions should be closed under
joins, which in the case of a categorical formulation leads to a
quantaloid structure. The inclusion 
$\q^-\a\hookrightarrow\q^{st}\a$ then follows by corollary
\ref{quantforMINCROSS}. 

Let us now apply all this to some particular categories that have
applications in physics. Consider  the subcategory $\j \c \l at_{\1}$ of
$\j \c \l at$ introduced in (Coecke and Moore, 1999) with as objects
complete lattices $L_i$  with `fixed' top $\1$ and an element
$1_i$ such that $\forall a_i\not=\1:a_i\leq 1_i$, and as morphisms join
preserving maps with $\1\mapsto\1$ or maps of which the image is
$\{0\}$. Since it is a subquantaloid of $\j \c \l at$, being itself also
a quantaloid for pointwise order, all the above considerations that
apply to
$\j \c \l at$ apply to $\j \c \l at_{\1}$.   
    
Let $\j \c \a \l at$ be the category of complete atomistic lattices with
as morphisms join-preserving maps that assign atoms to atoms or the
least element (Faure and Fr\"olicher, 1993; 1994). If we apply $\e$ to
this category we obtain a full subcategory $\e\j \c \a \l at$ of $\j \c
\l at$ as pre-enrichment. This is the minimal extension of $\j \c \a \l
at$ that assures that all
by $\q^-$ induced morphisms are state transitions. Since the category $\p
rop$ introduced in (Moore, 1995) is a full subcategory of $\j \c \a \l
at$ by restricting objects complete orthocomplemented atomistic lattices,
all the above applies to it, i.e., within this context one should rather
consider $\e\p rop$.
    
\section*{Acknowledgments} 

We thank F. Borceux, Cl.-A. Faure, D.J. Moore and J. Paseka for
discussing aspects of this paper. I.S. thanks VUB-FUND-DWIS for logistic
and financial support.

\section*{References}

\noindent Adamek, J., Herrlich, H. and Strecker, G.E. (1990) {\it
Abstract and Concrete Categories}, J. Wiley \& Sons Inc.  
   
\noindent Amira, H., Coecke, B. and Stubbe, I. (1998) {\it Helv. Phys.
Acta} {\bf 71}, 554.

\noindent Borceux, F. (1994) {\it Handbook of Categorical Algebra 1 and
2}, Cambridge University Press.  

\noindent Coecke, B. and Moore, D.J. (1999) `Operational Galois adjunctions', 
In\,: B. Coecke, D.J. Moore and A. Wilce, (Eds.), {\it Current Research in Operational
Quantum Logic: Algebras, Categories and Languages},    
pp.195--218, Kluwer Academic Publishers; quant-ph/0008021. 

\noindent Coecke, B. and Stubbe, I. (1999a) {\it Found. Phys. Lett.} {\bf
12}, 29; quant-ph/0008020.   
   
\noindent Coecke, B. and Stubbe, I. (1999b) {\it Int. J. Theor.
Phys.} {\bf 38}, 3269.
   
\noindent Faure, Cl.-A. and Fr\"olicher, A. (1993) {\it Geom. Dedicata}
{\bf 47}, 25.
 
\noindent Faure, Cl.-A. and Fr\"orlicher, A. (1994) {\it Geom. Dedicata}
{\bf 53}, 273.  
 
\noindent Faure, Cl.-A., Moore, D.J. and Piron, C. (1995) {\it Helv.
Phys. Acta} {\bf 68}, 150.  
   
\noindent Moore, D.J. (1995) {\it Helv. Phys. Acta} {\bf 68}, 658.
 
\noindent Moore, D.J. (1999) {\it Int. J. Theor.
Phys.} {\bf 38}, 793.

\noindent Piron, C. (1964) {\it Helv. Phys. Acta} {\bf 37}, 439.  

\noindent Piron, C. (1976) {\it Foundations of Quantum Physics}, W. A.
Benjamin, Inc.

\noindent Pitts, A.M. (1988) {\it Proc. London Math. Soc.} {\bf 57}, 433.

\noindent Pool, J.C.T. (1968) {\it Comm. Math. Phys.} {\bf 9}, 118.
 
\noindent Rosenthal, K.I. (1990)  {\it Quantales and their Applications},
Pitman Research Notes in Math. {\bf 234}, Longman.  

\noindent Rosenthal, K.I. (1991)  {\it J. Pure Appl. Alg.} {\bf 77},
67.  

\end{document}